\documentclass[11pt]{article}
\usepackage{fullpage}
\usepackage{amsmath}
\usepackage{amsthm}
\usepackage{graphicx}
\usepackage{pstricks}
\newpsobject{showgrid}{psgrid}{subgriddiv=1,griddots=10,gridlabels=6pt}
\newpsobject{showgrid}{psgrid}{subgriddiv=1,griddots=10,gridlabels=6pt}

\newcommand {\abs} [1] {\left| #1 \right|}

\newtheorem {thm} {Theorem}[section]

\numberwithin{equation}{section}

\title{Linear Stability Analysis of Resonant Periodic Motions in the Restricted Three-Body Problem}
\author{D. Viswanath \thanks{ Department of Mathematics, University of Michigan, 
530 Church Street, Ann Arbor, MI 48109, U.S.A.,
{\tt divakar@umich.edu}.
This work was supported by the N.S.F.
and by a research fellowship from the Sloan Foundation.}}

\begin{document}
\date{March 1, 2005}
\maketitle
\begin{abstract}
The equations of the restricted three-body problem describe the motion of a massless particle
under the influence of two primaries of masses $1-\mu$ and $\mu$, $0\leq \mu \leq 1/2$, that
circle each other with period equal to $2\pi$. When $\mu=0$, the problem admits  
orbits for the massless particle that are ellipses of eccentricity $e$
with the primary of mass $1$ located at one of the focii. If the period is a rational
multiple of $2\pi$, denoted $2\pi p/q$, some of these orbits perturb to periodic motions
for $\mu > 0$. For typical values of $e$ and $p/q$, two resonant periodic motions are obtained
for $\mu > 0$. We show that the characteristic multipliers of both these motions are given
by expressions of the form $1\pm\sqrt{C(e,p,q)\mu}+O(\mu)$ in the limit $\mu\rightarrow 0$.
The coefficient $C(e,p,q)$ is analytic in $e$ at $e=0$ and $C(e,p,q)=O(e^{\abs{p-q}})$. 
The coefficients in front of $e^{\abs{p-q}}$, obtained when $C(e,p,q)$ is expanded
in powers of $e$ for the two resonant periodic motions, sum to zero. Typically, if one of the two resonant periodic motions is of elliptic type the other
is of hyperbolic type. We give similar results for retrograde periodic motions and discuss
periodic motions that nearly collide with the primary of mass $1-\mu$. 
\end{abstract}

\section{Introduction}

Resonance between the periods of planets, satellites, and other objects is a feature of celestial mechanics. Among the many extra-solar planets discovered recently, many types of resonance have been observed. Here we consider resonant periodic motions in the
restricted three-body problem, with the principal aim of  gaining an understanding of the instabilities caused by resonance.

The restricted three-body problem is about the motion of a massless particle under the influence
of primaries of masses $1-\mu$ and $\mu$ that circle each other at a distance equal to $1$ and with period $2\pi$. 
We will assume $0\leq \mu \leq 1/2$ and in much of the discussion $\mu$ is small and positive.
Instead of the inertial frame of reference, it is often convenient to use a 
frame of reference which rotates with the two primaries and which is centered at the center of mass of the two primaries. The Hamiltonian for the motion of the massless particle
in a rotating Cartesian frame is 
\begin{equation}H = \frac{1}{2}(p_x^2+p_y^2)+y p_x - x p_y - \frac{1-\mu}{\Delta_0} -\frac{\mu}{\Delta_1},
\label{eqn-1-0}
\end{equation}
where
$\Delta_0 = ((x+\mu)^2+y^2)^{1/2}$ and $\Delta_1 = ((x-1+\mu)^2+y^2)^{1/2}$.
The generalized momenta corresponding to $x$ and $y$ are $p_x$ and $p_y$. The primaries of
mass $1-\mu$ and $\mu$ are located at $(-\mu,0)$ and $(1-\mu,0)$, respectively, and 
$\Delta_0$ and $\Delta_1$ are the distances from the massless particle 
to the two primaries. We consider the motion
only in regions where the distances $\Delta_0$ and $\Delta_1$ are both bounded away from
zero. Therefore $(1-\mu)/\Delta_0 + \mu/\Delta_1$ can be expanded in powers
of $\mu$ as 
\begin{equation}
\frac{1}{(x^2+y^2)^{1/2}} + \mu \Omega'(x,y,\mu)
            = \frac{1}{(x^2+y^2)^{1/2}} + \mu \Omega(x,y) + O(\mu^2)
\label{eqn-1-1}
\end{equation}
where 
\begin{equation*}
\Omega(x,y) = \frac{1}{((x-1)^2+y^2)^{1/2}}-\frac{x}{(x^2+y^2)^{3/2}}-
\frac{1}{(x^2+y^2)^{1/2}}.
\end{equation*}
The change of variables $p_x\leftarrow -p_x$, $x\leftarrow x$,
$p_y\leftarrow p_y$, $y\leftarrow -y$ is symplectic with multiplier $-1$ and leaves the Hamiltonian $H$
unchanged. Therefore if $x(t)=X(t), y(t)=Y(t)$, $0\leq t\leq T$, is a solution
of Hamilton's equations formed using $H$, then $x(t)=X(-t), y(t)=-Y(-t)$, $-T\leq t \leq 0$,
is also a solution. 

If $r$ and $\theta$ are polar coordinates then $x=r\cos\theta$, $y=r\sin\theta$, and
the Hamiltonian becomes 
\[H = \frac{1}{2}\Bigl(R^2 + \frac{G^2}{r^2}\Bigr)-G-\frac{1}{r}-\mu\Omega(r,\theta)+
O(\mu^2),\]
where
\begin{equation} 
\Omega(r,\theta) = \frac{1}{\Delta_1}-\frac{\cos\theta}{r^2}-\frac{1}{r}
=\frac{1}{(1+r^2-2r\cos\theta)^{1/2}}-\frac{\cos\theta}{r^2}-\frac{1}{r}.
\label{eqn-1-2}
\end{equation}
The discrete symmetry of this Hamiltonian is given by 
$R\leftarrow -R$, $r\leftarrow r$, $G\leftarrow G$, and $\theta\leftarrow -\theta$. This transformation
too is symplectic with multiplier $-1$.

Given $R,G,r,\theta$, another set of variables better suited for perturbative calculations
can be defined as follows. Let $H_0 = (R^2+G^2/r^2)/2-1/r$. Define the variables
$a,e,E,\nu,g$ using the equations
\begin{align}
e &= \sqrt{1+2G^2H_0}\nonumber\\
a(1-e) &= G^2/(1+e) \nonumber\\
r &= a(1-e \cos E) \nonumber\\
\cos \nu &= \frac{\cos E -e}{1-e\cos E}\nonumber\\
g &= \theta - \nu. \label{eqn-1-3}
\end{align}
When the third and fourth equations in \eqref{eqn-1-3} are solved for $E$ and $\nu$, there
are typically two solutions in the interval $[0,2\pi)$. Exactly one of these can be chosen
by heeding the signs of $G$ and $R$. Until the last section, we assume $0<e<1$.  
The variables $a,e,\nu,g$ can be used instead of $R,G,r,\theta$. When $\mu=0$, the orbit
(of the massless particle) is an ellipse in the inertial frame with eccentricity $e$ and
semimajor axis $a$; the angles $E$ and $\nu$ are the eccentric and true anomalies, 
respectively; and $g$ is the argument of the perihelion in the rotating frame.

The change of variables from $R,G,r,\theta$ to $a,e,\nu,g$, in the region of phase space
with $0<e<1$, is not symplectic. To obtain a symplectic change of variables in that
region, define the variables $L,l$ using the equations
\begin{equation}
L = \pm\sqrt{a}\quad\quad l = E-e\sin E.
\label{eqn-1-4}
\end{equation}
As further explained in Section 3, the sign of $L$ is chosen to be the same as that of
$G$, the angular momentum. The change to $L,G,l,g$ is symplectic and the Hamiltonian
becomes
\begin{equation}
 H = -\frac{1}{2L^2} - G - \mu\Omega(L,G,l,g) + O(\mu^2).
 \label{eqn-1-5}
\end{equation}
Simple expressions for the perturbation $\Omega$ in terms of $x,y$ and in terms of
$r,\theta$ have been given in \eqref{eqn-1-1} and \eqref{eqn-1-2}.
In Sections 2 and 3, we partially develop $\Omega$ as a 
trigonometric series in $l$ and $g$. The variable $l$ is called the mean anomaly.
The discrete symmetry of the Hamiltonian \eqref{eqn-1-5} is given by
$L\leftarrow L, l\leftarrow -l, G\leftarrow G, g\leftarrow -g$. This transformation too
is symplectic with multiplier $-1$.

\begin{figure}
\begin{center}
\begin{pspicture}(-3,-3)(3,3)
\psline[arrowsize=3pt 2]{->}(0,-3)(0,3)
\psline[arrowsize=3pt 2]{->}(-3,0)(3,0)
\rput{-20}(-0.23,.083){\psellipse(0,0)(2,1.5) \pscircle[fillstyle=solid,fillcolor=white](1,1.3){0.1}
          \pscurve[arrowsize=6pt 2]{->}(1,1.3)(0.8,1.37)(0,1.5)}
\psarcn[arrowsize=6pt 2]{->}(0,0){2.0}{10}{300}
\pscircle[linestyle=dotted]{2.6}
\pscircle*(2.6,0){0.1}
\pscircle*{0.15}
\uput[u](2.9,0){$x$}
\uput[r](0,2.9){$y$}
\end{pspicture}
\end{center}
\caption{The dotted line above is a circle of radius $1$.}
\label{fig-1}
\end{figure}
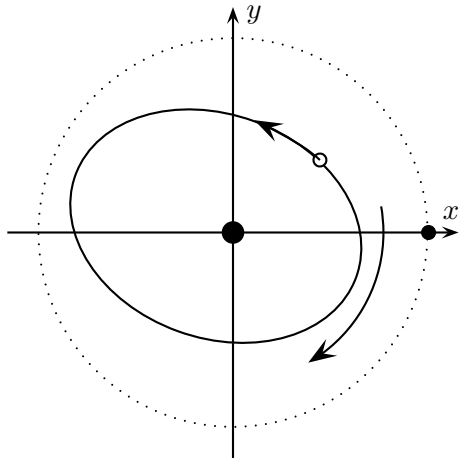

Let $\mu=0$.
If $L=(p/q)^{1/3}$, where $p$ and $q$ are relatively prime positive integers, or equivalently
$a=(p/q)^{2/3}$, and $0<e<1$ or $G=(a(1-e^2))^{1/2}$, then the massless particle moves on
an ellipse in the counterclockwise sense. In the rotating frame, this ellipse itself rotates
in the clockwise sense with angular velocity equal to $1$ as shown in Figure \ref{fig-1}. 
The orbit has a period equal to $2\pi p$ in the rotating frame as the period of motion on the
ellipse is in resonance with the period of rotation of the frame. 
Arenstorf \cite{Arenstorf1}
and Barrar \cite{Barrar1} proved that the periodic orbit persists for $\mu > 0$ and
$\mu$ small enough, if $l$ and $g$ are multiples of $\pi$ at $t=0$, and if the orbit for $\mu=0$ does not pass through the location of the primary
of mass $\mu$. 

\begin{figure}
\begin{center}
\begin{pspicture}(-3,-3)(3,3)
\psline[arrowsize=3pt 2]{->}(0,0)(3,0)
\psline(0,0)(2.12,2.12)
\psarc{->}{0.5}{0}{45}
\pscircle[linestyle=dashed](0,0){2.2}
\pscircle[linestyle=dashed](0,0){2.7}
\uput[ur](0.6,0){$l$}
\uput[u](2.9,0){$L$}
\qline(2.25,0.2)(2.65,-0.2)
\qline(2.25,-0.2)(2.65,0.2)
\pscircle(-2.45,0){0.12}
\qdisk(-2.45,0){0.03}
\end{pspicture}
\hspace{1.5cm}
\begin{pspicture}(-3,-3)(3,3)
\psline[arrowsize=3pt 2]{->}(0,0)(3,0)
\psline(0,0)(2.12,2.12)
\psarc{->}{0.5}{0}{45}
\pscircle[linestyle=dashed](0,0){2.2}
\pscircle[linestyle=dashed](0,0){2.7}
\uput[ur](0.6,0){$l$}
\uput[u](2.9,0){$L$}

\qdisk(2.45,0){0.03}
\pscircle(2.45,0){0.12}
\qdisk(-2.45,0){0.03}
\pscircle(-2.45,0){0.12}
%\qline(2.25,0.2)(2.65,-0.2)
%\qline(2.25,-0.2)(2.65,0.2)
%\qline(-2.25,0.2)(-2.65,-0.2)
%\qline(-2.25,-0.2)(-2.65,0.2)

\qline(0.2,2.25)(-0.2,2.65)
\qline(-0.2,2.25)(0.2,2.65)
\qline(0.2,-2.25)(-0.2,-2.65)
\qline(-0.2,-2.25)(0.2,-2.65)
%\qdisk(0,2.45){0.03}
%\pscircle(0,2.45){0.12}
%\qdisk(0,-2.45){0.03}
%\pscircle(0,-2.45){0.12}

\end{pspicture}
\end{center}
\caption{Periodic points on the $L$-$l$ plane with $p=1$ and $p=2$, respectively.}
\label{fig-2}
\end{figure}
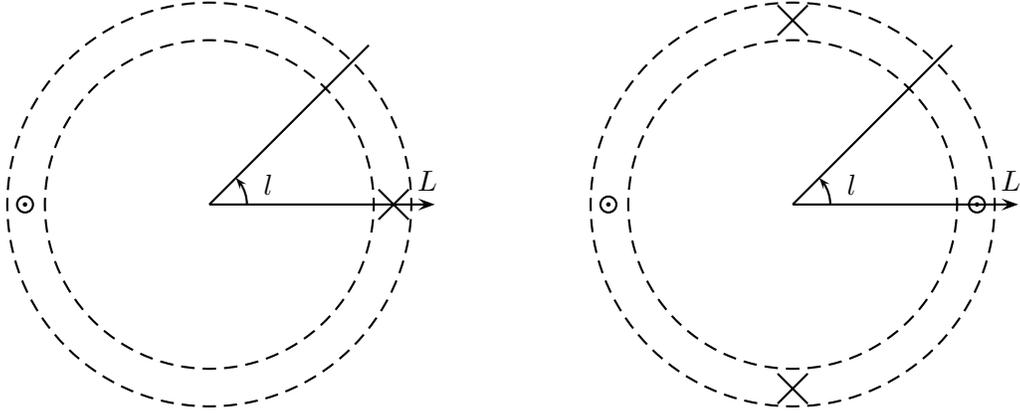

Since $l$ and $g$ are angular variables, there are four possibilities for $l$ and $g$ at
$t=0$, namely $l=n_l \pi$ and $g=n_g \pi$ with $n_l$ and $n_g$ either $0$ or $1$. Only two
of these are distinct. They are depicted on the Poincar\'e section $g=0$ in Figure \ref{fig-2}.
If $p$ is odd the two distinct orbits are obtained with $n_g=0$ and $n_l$ either $0$ or $1$.
When $p$ is even, we can take $n_l=0$ and $n_g$ either $0$ or $1$. Over its period, the orbit
intersects the Poincar\'e section exactly $p$ times. The plots in Figure \ref{fig-2}
correspond to $p=1$ and $p=2$. 

We investigate the linear stability of resonant periodic motions, such as the ones
depicted in Figure \ref{fig-2}, for $\mu > 0$. The existence and stability of the nearly circular orbits, which are related to the nearly elliptic ones considered here, is discussed in
\cite{Birkhoff1}, \cite{LS1}, \cite{MH1}, and \cite{SM1}.

In Section 2, we prove that the characteristic multipliers of the resonant periodic orbits
are given by the expression $1\pm\sqrt{C(e,p,q)\mu}+O(\mu)$ in the limit
$\mu\rightarrow 0$ where
\[C(e,p,q) = -\frac{6\pi p^{5/3}}{q^{2/3}}\int_0^{2\pi p} 
\frac{\partial^2 \Omega(r,\theta)}{\partial \theta^2} dt.\]
In the integral above, expression \eqref{eqn-1-2}  for $\Omega(r,\theta)$
must be used to evaluate the partial derivatives of $\Omega(r,\theta)$
with respect to $\theta$. After that, the variables $r$ and $\theta$ must be expressed in terms of
$L,G,l,g$ using \eqref{eqn-1-3} and \eqref{eqn-1-4}. The dependence of $l$ and $g$
on $t$ is obtained from \eqref{eqn-1-5} with $\mu=0$ as $l(t) = qt/p +n_l\pi$ and
$g(t) = -t + n_g\pi$, and $L=(p/q)^{1/3}$ and $G=(p/q)^{1/3}(1-e^2)^{1/2}$ do not vary
with $t$.
 
If $p/q\neq 1/1$, $C(e,p,q)$ is analytic at $e=0$. In Section 2, we prove that all
powers of $e$ with exponent less than $\abs{p-q}$ have coefficients equal to $0$
in the power series expansion of $C(e,p,q)$ about $e=0$. The coefficient of $e^{\abs{p-q}}$ can be obtained as
a finite sum of hypergeometric terms if $q>1$. The coefficients of $e^{\abs{p-q}}$ for the two
distinct types determined by the choice of $n_l$ and $n_g$ sum to zero. It seems to be 
the typical case that if one of the two types of periodic motions is elliptic  then
the other is hyperbolic for small $\mu$. This is depicted in Figure \ref{fig-2} 
by the use of crosses
for hyperbolic points and circles for elliptic points.  In later work, we will 
demonstrate 
the existence of homoclinic points in the resonance band around $L=(p/q)^{1/3}$.

In Section $3$, we give similar results for retrograde periodic motions. 
For retrograde periodic motions, $C(e,p,q) = O(e^{p+q})$ as $e\rightarrow 0$.
In Section 4,
we consider periodic motions near collision with the primary of mass $1-\mu$. The 
variables $L,G,l,g$ are not a valid choice of variables at collision. Giacaglia
\cite{Giacaglia1} solved the Hamilton-Jacobi equation for the regularized Hamiltonian
with $\mu=0$
and derived action-angle variables that are valid near collision. The existence of 
resonant periodic motions near collision was proved by Schmidt \cite{Schmidt1}. 
We solve the Hamilton-Jacobi equation again, using the geometric approach of 
Arnold \cite{Arnold1}, to correct a calculation and to put it on a more secure
basis. We discuss the existence of periodic motions that  collide with the primary of mass $1-\mu$.

\section{Direct periodic motions}
Let $\mu = 0$. 
The Hamilton's equations formed using \eqref{eqn-1-5} are
$\dot{l} = L^{-3}$, $\dot{g} = -1$, $\dot{L} = 0$, and $\dot{G} = 0$.
Consider the solution with the initial conditions $L(0)=(p/q)^{1/3}$,
$G(0)=(p/q)^{1/3}(1-e^2)^{1/2}$, $l(0)=n_l \pi$, and $g(0)=n_g\pi$, where $p$ and $q$ are relatively prime positive integers, $0<e<1$, and $n_l,n_g$ are either $0$ or $1$. This solution has a period equal to $2\pi p$. Assume that $\Delta_1 > 0$ along this orbit. Then $\Omega'$ in
\eqref{eqn-1-1} is uniformly bounded in a neighborhood of the orbit for $\mu$ small enough. Therefore, for initial conditions chosen in a neighborhood of this solution and for a finite interval of time, the solution is analytic in $\mu$ at $\mu=0$. The assumption
$\Delta_1 > 0$ is valid for $p/q \neq 1/1$ if $e$ is small enough. For given $p/q$, $n_l$, and $n_g$, it is valid for all except finitely many values of $e$ \cite{Arenstorf1}.

The solution is unchanged by the discrete symmetry $L\leftarrow L$,
$l\leftarrow-l$, $G\leftarrow G$, $g\leftarrow-g$, at both $t=0$ and $t=\pi p$.
For $\mu > 0$ and $\mu$ sufficiently small, it is possible to perturb $L(0)$ and the time of flight $\pi p$ such that the solution is unchanged by the discrete symmetry at 
$t=0$ and $t=\pi p + O(\mu)$ \cite{Barrar1}.  The existence of a unique family of periodic solutions
for $\mu > 0$ and $\mu$ sufficiently small
with period $T = 2\pi p + O(\mu)$ and initial conditions $L(0)=(p/q)^{1/3}+O(\mu)$, $G(0) = (p/q)^{1/3}(1-e^2)^{1/2}$, $l(0) = n_l\pi$, and $g(0) = n_g\pi$ follows.

\subsection{Characteristic multipliers}

To obtain the characteristic multipliers of this family of periodic solutions, we consider the return map to the Poincar\'e section $g = n_g\pi$. 
The Hamilton's equations
that correspond to \eqref{eqn-1-5} are 
\begin{equation}
\dot{l} = L^{-3}-\mu \Omega_L + O(\mu^2),\quad
\dot{g} = -1-\mu\Omega_G + O(\mu^2),\quad
\dot{L} = \mu \Omega_l+O(\mu^2),\quad
\dot{G} = \mu \Omega_g + O(\mu^2).
\label{eqn-2-1}
\end{equation}
Since $H$ is a first integral of these equations, we can fix $H$ and solve \eqref{eqn-1-5} for $G$ using the implicit function theorem and identify the Poincar\'e section
with the $L$-$l$ plane.
Figure \ref{fig-2} depicts the intersections of periodic solutions with $p=1$
and $p=2$ with the Poincar\'e section.  We refer to these intersections as periodic
points.
Since $\dot{g} = -1$ for $\mu=0$, the implicit function theorem implies the existence of the return map over large regions of the $L$-$l$ plane and certainly in a neighborhood
of the periodic points. The Poincar\'e map leaves the area element $dl\, dL$ invariant. 
When $\mu=0$, the return time is $2\pi$ for all points in the $L$-$l$ plane.

When $\mu=0$, the $p$th return map maps the  periodic point on the $L$-$l$ section
with $L(0)=(p/q)^{1/3}$ and $l(0) = n_l\pi$ to itself. The $p$th return time
denoted by $T$ is $2\pi p$. By \eqref{eqn-2-1}, $l(t) = l(0) + t L(0)^{-3}$
and $L(t) = L(0)$. Therefore,
\begin{equation}
\frac{\partial(L(T), l(T))}{\partial(L(0), l(0)} =
\begin{pmatrix} 1 & 0\\ -\frac{6\pi q^{4/3}}{p^{1/3}} & 1 \end{pmatrix}
\label{eqn-2-2}
\end{equation}
at $L(0)=(p/q)^{1/3}$, $l(0) = n_l\pi$, which implies that the characteristic
multipliers of the $p$th return map are both $1$. Two of the characteristic 
multipliers of any periodic solution of \eqref{eqn-2-1} must be $1$. In this case
all four characteristic multipliers equal $1$.

The periodic solution with the initial conditions $L(0)=(p/q)^{1/3}+O(\mu)$, $G(0) = (p/q)^{1/3}(1-e^2)^{1/2}$, $l(0) = n_l\pi$, and $g(0) = n_g\pi$ depends analytically
on $\mu$ at $\mu=0$ and the $p$th return time, which is also the period,
is given by $T=2\pi p + O(\mu)$. Therefore the
entries of the Jacobian matrix in \eqref{eqn-2-2} will be perturbed by $O(\mu)$ for
$\mu > 0$ and $\mu$ sufficiently small. We determine only the perturbation to the entry
in the upper right corner which is equal to $0$ in \eqref{eqn-2-2}.  
By \eqref{eqn-2-1}, 
\[L(T) = L(0) + \mu \int_{0}^T \Omega_l(L(t),G(t),l(t), g(t))dt
+ O(\mu^2).\]
As the solution depends analytically on $\mu$ at $\mu=0$, we have 
\begin{align*}L(t) &= (p/q)^{1/3}+ O(\mu)\\
 G(t) &= (p/q)^{1/3}(1-e^2)^{1/2} + O(\mu)\\ 
l(t) &= n_l\pi + qt/p+O(\mu)\\ 
g(t)&= n_g\pi-t+O(\mu)
\end{align*} 
for $0\leq t\leq T$.
The $p$th return time $T$ itself can depend upon $l(0)$. However,
$\frac{\partial T}{\partial l(0)}=0$ at $\mu=0$ as the $p$th return time equals
$2\pi p$ for all points on the $L$-$l$ plane when $\mu=0$. Therefore
\[\frac{\partial L(T)}{\partial l(0)} = \mu \int_0^{2\pi p} \Omega_{ll}\bigl(
(p/q)^{1/3},\, (p/q)^{1/3}(1-e^2)^{1/2},\, n_l\pi+qt/p,\,n_g\pi-t\bigr) dt + O(\mu^2).\]
Define $C(e,p,q)$ such that the integral in the above equation is given
by $-C(e,p,q) p^{1/3}/(6\pi q^{4/3})$.
Then the Jacobian matrix for the $p$th return map is given by
\begin{equation*}
\frac{\partial(L(T), l(T))}{\partial(L(0), l(0)} =
\begin{pmatrix} 1 & 0\\ -\frac{6\pi q^{4/3}}{p^{1/3}} & 1 \end{pmatrix}
+\mu \begin{pmatrix} \star & -C(e,p,q)\frac{p^{1/3}}{6\pi q^{4/3}}\\
\star & \star
\end{pmatrix}
+O(\mu^2).
\end{equation*}
It follows that the characteristic multipliers are $1\pm \sqrt{C(e,p,q)\mu}
+ O(\mu)$ in the limit $\mu\rightarrow 0$. If $C(e,p,q)>0$, the periodic solutions
in the family are hyperbolic for $\mu>0$ and $\mu$ sufficiently small. 
If $C(e,p,q)<0$, those periodic solutions are elliptic.

The expression for $C(e,p,q)$, namely
\begin{equation}
C(e,p,q) = -\frac{6\pi q^{4/3}}{p^{1/3}} \int_0^{2\pi p}
\Omega_{ll}\bigl((p/q)^{1/3},\,(p/q)^{1/3}(1-e^2)^{1/2},\, n_l \pi + qt/p,\, n_g\pi-t
\bigr) dt
\label{eqn-2-3}
\end{equation}
can be cast into a simpler form. If $\Omega$ has the same arguments as in 
\eqref{eqn-2-3}, it is periodic in $t$ with period $2\pi p$. Using
the identity $\frac{d\Omega_l}{dt} = \frac{q}{p} \Omega_{ll} - \Omega_{lg}$,
the integrand in \eqref{eqn-2-3} can be changed to $\Omega_{lg}$. Use of the
identity $\frac{d\Omega_g}{dt} = \frac{q}{p}\Omega_{lg}-\Omega_{gg}$ gives
the expression
\begin{equation}
C(e,p,q) = -\frac{6\pi p^{5/3}}{q^{2/3}}\int_0^{2\pi p}
\Omega_{gg}\bigl((p/q)^{1/3},\,(p/q)^{1/3}(1-e^2)^{1/2},\, n_l \pi + qt/p,\, n_g\pi-t
\bigr) dt.
\label{eqn-2-4}
\end{equation}
The advantage of \eqref{eqn-2-4} becomes evident when the expression 
\eqref{eqn-1-2} for $\Omega$ in terms of $r$ and $\theta$ is considered.
From \eqref{eqn-1-3} and \eqref{eqn-1-4}, we see that $r$ depends on
$L$, $G$, and $l$ but not on $g$ and that $\theta = \nu - g$, where
$\nu$ too depends upon $L$, $G$, and $l$ but not on $g$. Therefore
\begin{equation}
C(e,p,q) = -\frac{6\pi p^{5/3}}{q^{2/3}}\int_0^{2\pi p}
\Bigl(\frac{1}{\Delta_1}\Bigr)_{\theta\theta} + \frac{\cos \theta}{r^2} dt,
\label{eqn-2-5}
\end{equation}
where $\Delta_1 = (1+r^2-2r\cos\theta)^{1/2}$. The variables $r$ and $\theta$
can be obtained in terms of $L,G,l,g$ using \eqref{eqn-1-3} and \eqref{eqn-1-4}.
For the family of periodic solutions determined by $p$, $q$, $e$, $n_l$, and
$n_g$, the variables $L,G,l,g$ depend upon $t$ as indicated by the arguments
of $\Omega_{ll}$ and $\Omega_{gg}$ in \eqref{eqn-2-3} and \eqref{eqn-2-4},
respectively.

For the family of periodic solutions to exist, we required $e\in (0,1)$ to be such
that $\Delta_1>0$ everywhere along the unperturbed orbit at $\mu=0$. The integrand of \eqref{eqn-2-5} can be differentiated with
respect to $e$ in a complex neighborhood of the value of $e$
for $0\leq t \leq 2\pi p$ using formulas, since $\Delta_1 > 0$, and is therefore analytic in
$e$. That $C(e, p, q)$ is also analytic in $e$ follows
from a standard argument (see \cite{Goursat1}).
If
$p/q \neq 1/1$, then $\Delta_1 > 0$ along the  unperturbed orbit for $e$
sufficiently small. Therefore $C(e,p,q)$ is analytic  at $e=0$ if $p/q \neq 1/1$.
The expression \eqref{eqn-2-5} can be used to investigate the expansion of
$C(e,p,q)$ in powers of $e$ assuming $p/q\neq 1/1$.

\subsection{Expansion of $C(e,p,q)$}

Expression \eqref{eqn-2-5} for $C(e,p,q)$ is an
integral over an elliptic solution of the two-body problem,
obtained by setting $\mu=0$ in \eqref{eqn-2-1},
that satisfies $\Delta_1 > 0$ at all points on the solution. By \eqref{eqn-2-1}
and \eqref{eqn-1-4}, $l = qt/p+n_l\pi$ and $l = E-e\sin E$. Change the
variable of integration in \eqref{eqn-2-5} to $E$ using $dt = (p/q)^{1/3}r\,dE$,
and then to $F$, which is defined by $qF=E$,
to get
\begin{equation}
C(e,p,q) = -\frac{6\pi p^{2}}{q}\int_0^{2\pi q} 
\Bigl(\frac{r}{\Delta_1}\Bigr)_{\theta\theta}+\frac{\cos\theta}{r}\, dE
= -6\pi p^{2} \int_0^{2\pi} \Bigl(\frac{r}{\Delta_1}\Bigr)_{\theta\theta}
+\frac{\cos\theta}{r}\, dF.
\label{eqn-2-6}
\end{equation}
Define $C_1(e,p,q)$ and $C_2(e,p,q)$ by
\begin{equation}
C_1(e,p,q) = \int_0^{2\pi}\Bigl(\frac{r}{\Delta_1}\Bigr)_{\theta\theta}\, dF\quad
\text{and}\quad C_2(e,p,q) = \int_0^{2\pi} \frac{\cos\theta}{r}\, dF.
\label{eqn-2-7}
\end{equation}
Then $C_(e,p,q) = -6\pi p^{2}(C_1(e,p,q)+C_2(e,p,q))$
The dependence of $r$ and $\theta$ on $E$ or $F$ is given by
$r=(p/q)^{2/3}(1-e\cos E)$ and $\theta=\nu + n_g\pi -t$. The variables
$t$ and $\nu$ can be obtained as functions of $E$ or $F$ using
$l = qt/p+n_l\pi$, \eqref{eqn-1-3}, and \eqref{eqn-1-4}.
The integrands of both integrals in \eqref{eqn-2-7} are periodic in $F$
with period $2\pi$, and as $\Delta_1 > 0$, are analytic in $F$ for
$0\leq F < 2\pi$.

Let $z=\exp(iF)$. Define $\beta$ by $e = 2\beta/(1+\beta^2)$.
Standard formulas that connect the true and mean anomalies
with the eccentric anomaly (see \cite{BC1}) imply the following:
\begin{align}
\exp{i\nu} &= z^q (1-\beta z^{-q})(1-\beta z^q)^{-1}\nonumber\\
\exp{it} &= \exp(-i\pi n_lp/q)
z^p \exp\Bigl(\frac{ep}{2q}(z^{-q}-z^q)\Bigl).
\label{eqn-2-8}
\end{align}
The variable $r$ and $\exp(in\theta)$, $n\in Z$, too can be expressed in terms of $z$.
\begin{align}
r &= (p/q)^{2/3}(1+\beta^2)^{-1}(1-\beta z^{-q})(1-\beta z^{q})\nonumber\\ 
\exp(in\theta) &= 
(-1)^{n n_g}\exp(i\pi n n_lp/q)
z^{n(q-p)}(1-\beta z^{-q})^n (1-\beta z^q)^{-n}
\exp\Bigl(\frac{enp}{2q}(z^q-z^{-q})\Bigr).
\label{eqn-2-9}
\end{align}
If the integrands in \eqref{eqn-2-7} are expressed in terms of $z$, they are analytic
on the circle $\abs{z} = 1$. Therefore the integrands can be expanded in  Laurent
series in $z$. Only the constant terms in the Laurent series contribute
to $C_1(e,p,q)$ and $C_2(e,p,q)$. For later use, we record the following Laurent
series:
\begin{equation}
\exp\Bigl(\frac{enp}{2q}(z^q-z^{-q})\Bigr)
= \sum_{k=-\infty}^{\infty} J_k\bigl(enp/q\bigr) z^{kq}.
\label{eqn-2-10}
\end{equation}
The Bessel functions $J_k$ are entire functions. We need the first term in
the power series expansion of $J_k(x)$ about $x=0$, which is $\frac{x^k}{2^k k!}$
for $k\geq 0$,
but no others. For $k>0$, $J_{-k}(x) = (-1)^k J_k(x)$.

When we consider $C_2(e,p,q)$, defined by \eqref{eqn-2-7},
and use a Laurent series of its integrand in $z$,
we find that $C_2(e, p, q)\equiv 0$ except when $q=1$. Since $\cos\theta/r
= (\exp(i\theta)+\exp(-i\theta))/(2r)$, we may use expressions for
$\exp(in\theta)$, $n=\pm 1$, and $r$ from \eqref{eqn-2-9} to get
\[\frac{\cos\theta}{r} = \Bigl(\frac{q}{p}\Bigr)^{2/3} \frac{1+\beta^2}{2}
(\exp(i\pi n_lp/q)z^{q-p} F(z) + \exp(-i\pi n_lp/q)z^{p-q}F(1/z)),\]
where \(\quad F(z) = (-1)^{n_g}\,(1-\beta z^q)^{-2} \exp\bigl(((ep)/(2q)) (z^q-z^{-q})\bigr).\)
The Laurent series of $F(z)$ only has powers of $z^q$. Since $p$ and $q$ are
relatively prime, none of these powers can equal $z^{p-q}$ or $z^{q-p}$,
except when $q=1$. Therefore the Laurent series of $\cos\theta/r$ has no constant
term if $q\neq 1$ and $C_2(e,p,q)=0$.

If $q=1$, we can use the Laurent series \eqref{eqn-2-10} with $n=1$ and
the binomial series of $(1-\beta z^{q})^{-2}$ to obtain the Laurent series
of $F(z)$ and of $\cos\theta/r$. The Laurent series of $\cos\theta/r$ implies
\begin{equation}
C_2(e,p,q) = (-1)^{n_g}(-1)^{n_l p}\,2\pi\frac{1+\beta^2}{p^{2/3}} \bigl(J_{p-1}(ep)
+2J_{p-2}(ep)\beta + 3J_{p-3}(ep)\beta^2+\cdots\bigr)
\label{eqn-2-11}
\end{equation}
for $q=1$. Note that $C_2(e,p,q) = O(e^{p-1})$ as $e\rightarrow 0$, when
$q=1$.

We will prove that the power series of $C_1(e, p, q)$, defined by \eqref{eqn-2-7},
has no terms lower than $e^{\abs{p-q}}$. We will also determine the 
coefficient of the $e^{\abs{p-q}}$ term. In the ensuing analysis, we assume
$p<q$. The $p>q$ case is similar. We obtain terms in the Laurent series of 
the integrand $\bigl(\frac{r}{\Delta_1}\bigr)_{\theta\theta}$ in two steps.
For the first step, consider
\[I(\alpha) = \frac{\partial^2}{\partial\theta^2}\Bigl(\frac{\alpha}
{\sqrt{1+\alpha^2-2\alpha\cos\theta}}\Bigr).\]
This quantity $I(\alpha)$ would equal the integrand if $\alpha = r$. However, we assume
$\alpha$ to be fixed in the range $(0,1)$. It will be set equal to $(p/q)^{2/3}$
later. The Fourier expansion of $I(\alpha)$ is 
\begin{equation}
I(\alpha) = \frac{1}{2} \sum_{n=-\infty}^{\infty} -n^2\alpha b_n(\alpha) \exp(in\theta),
\label{eqn-2-12}
\end{equation}
where $b_n(\alpha)$ are the Laplace coefficients. The Laplace coefficients are 
hypergeometric functions whose series converge for $\abs{\alpha} < 1$. They
and their derivatives satisfy a number of identities. See \cite{BC1}. 

In the second step, we let $\alpha = r$ in $I(\alpha)$ to make it equal to
the integrand $\bigl(\frac{r}{\Delta_1}\bigr)_{\theta\theta}$ in 
\eqref{eqn-2-7}. Let $f(x)$ be analytic in $x$ at $x=x_0$. Then $f((1+\delta)x_0)$
has a convergent power series expansion in $\delta$.  This power series
can be conveniently represented as $(1+\delta)^D f(x_0)$, with the understanding
that $D$ stands for the differential operator $x\frac{d}{dx}$ and that 
$(1+\delta)^D$ is expanded according to the binomial formula \cite{BC1}. 
Note that $r = (p/q)^{2/3}(1-e\cos E)$. Let $D$ be the differential operator
$\alpha \frac{d}{d\alpha}$. Since $\cos E = (z^{-q}+z^{q})/2$, we get
\begin{equation}
D_\alpha = (1-e\cos E)^D = (1+\beta^2)^{-D}(1-\beta z^{-q})^D(1-\beta z^q)^D,
\label{eqn-2-13}
\end{equation}
where the first equality defines the operator $D_\alpha$. 
By \eqref{eqn-2-12} and \eqref{eqn-2-13},
the integrand 
in the definition \eqref{eqn-2-7} of $C_1(e,p,q)$  can be represented as
\begin{equation}
\Bigl(\frac{r}{\Delta_1}\Bigr)_{\theta\theta}
= \frac{1}{2}\sum_{n=-\infty}^{\infty} -n^2 D_\alpha \alpha b_n(\alpha)\exp(in\theta),
\label{eqn-2-14}
\end{equation}
which is to be evaluated at $\alpha = (p/q)^{2/3}$. We will use \eqref{eqn-2-9} and 
\eqref{eqn-2-13} to substitute expressions in terms of $z$ for
$\exp(in\theta)$ and $D_\alpha$, respectively.

Define $X_n(A,B,C)= (1+\beta^2)^A(1-\beta z^{-q})^B (1-\beta z^q)^C
\exp\bigl((enp/(2q))(z^q-z^{-q})\bigr)$. The formal Laurent series of
$X_n(A,B,C)$ is gotten by multiplying the binomial expansions of
$(1-\beta z^{-q})^B$ and $(1-\beta z^q)^C$ with the series in \eqref{eqn-2-10}.
Each of these three series is in terms of integral powers of $z^q$ and so is the
resulting series for $X_n(A, B, C)$. In each of the three series the coefficient
of $z^{kq}$ is a power series in $e$ with the lowest power of $e$ with 
possibly nonzero
coefficient being $e^{\abs{k}}$ and so it is for the series expansion
of $X_n(A,B,C)$.

The $n$th term in the summation
in \eqref{eqn-2-14} is given by
\begin{equation}
-n^2 D_\alpha \alpha b_n(\alpha)\exp(in\theta) =
-(-1)^{nn_g} \exp(i\pi nn_lp/q) n^2 z^{n(q-p)}X_n(-D,D+n,D-n)(\alpha b_n(\alpha)).
\label{eqn-2-15}
\end{equation} 
For the Laurent series of the 
term in \eqref{eqn-2-15} above to have a nonzero constant term, $n$ must be an
integer multiple
of $q$ with $n\neq 0$. Let $n=n_0 q$. Then that constant term is $\pm n^2$ times the 
coefficient of $z^{n_0 (p-q) q}$ in the Laurent series of
$X_n(-D,D+n, D-n)$. By the previous paragraph, the lowest term 
with possibly nonzero coefficient in the power
series expansion of that constant term is $e^{\abs{n_0(p-q)}}$. As
$C_1(e,p,q)$ is obtained in \eqref{eqn-2-7} by integrating the sum
in \eqref{eqn-2-14}, only the constant term in the Laurent series in $z$
of each term of \eqref{eqn-2-14} makes a contribution to $C_1(e,p,q)$.
Thus the lowest possibly nonzero term in the power series of $C(e,p,q)$ is
$e^{\abs{p-q}}$ and the only terms of \eqref{eqn-2-14} that contribute to the
coefficient of $e^{\abs{p-q}}$ are $n=q$ and $n=-q$.

Thus we have proved that the coefficients of all powers of $e$ less than
$\abs{p-q}$ in the power series of $C_1(e,p, q)$ about $e=0$ are zero,
if $p<q$. The coefficient of $e^{\abs{p-q}}$ is equal to
\begin{equation}
-2\pi q^2 \frac{(-1)^{q-p}}{2^{q-p}}(-1)^{n_g q+n_l p}
\Biggl(\sum_{k=0}^{q-p}\binom{D+q}{k}
\frac{p^{q-p-k}}{(q-p-k)!}\Biggr)(\alpha b_q(\alpha))
\label{eqn-2-16}
\end{equation}
evaluated at $\alpha = (p/q)^{2/3}$.  By  definition $D$ is the
operator $\alpha \frac{d}{d\alpha}$, and $b_q$ is a hypergeometric function
whose series converges for $\abs{\alpha} < 1$ as mentioned earlier.

If $p > q$, a similar analysis proves that all powers of $e$ less than
$\abs{p-q}$ in the power series of $C_1(e,p, q)$ about $e=0$ are zero. 
The coefficient of $e^{\abs{p-q}}$  in that power series 
is equal to
\begin{equation}
-2\pi q^2 \frac{(-1)^{p-q}}{2^{p-q}}(-1)^{n_g q+n_l p}
\Biggl(\sum_{k=0}^{p-q}(-1)^k\binom{-D-q}{k}
\frac{p^{p-q-k}}{(p-q-k)!}\Biggr)b_q(\alpha),
\label{eqn-2-17}
\end{equation}
evaluated at $\alpha = (q/p)^{2/3}$.
That $C(e,p,q) = O(e^{\abs{p-q}})$ for
$p/q\neq 1/1$ is now clear from \eqref{eqn-2-6},
\eqref{eqn-2-7}, \eqref{eqn-2-11}, and the power series of $C_1(e,p,q)$ given
by \eqref{eqn-2-16} and \eqref{eqn-2-17}.

The two families of periodic solutions of \eqref{eqn-2-1} with
$L(0) = (p/q)^{1/3}+O(\mu)$ and $G(0) = (p/q)^{1/3}(1-e^2)^{1/2}+O(\mu)$
are given by $l(0) = n_l \pi$ and $g(0) = n_g \pi$, with $n_g=0$
and $n_l$ either $0$ or $1$ if $p$ is odd, and with $n_l=0$
and $n_g$ either $0$ or $1$ if $p$ is even. Since $p$ and $q$ are
relatively prime, $q$ must be odd if $p$ is even. By \eqref{eqn-2-11},
\eqref{eqn-2-16}, and \eqref{eqn-2-17}, the coefficients of $e^{\abs{p-q}}$
in the power series of $C(e,p,q)$ for the two families sum to zero.
Of course, there is still the possibility that both these coefficients are
zero. This possibility can be eliminated by evaluating \eqref{eqn-2-16} if
$p<q$ and \eqref{eqn-2-17} if $p>q$, and also \eqref{eqn-2-11} if $q=1$.

\subsection{Statement of results}

The theorem below summarizes the results of this section.

\begin{thm}
There exists a $\mu$-dependent family of periodic solutions of 
Hamilton's equations
\eqref{eqn-2-1} with the initial conditions $L(0) = (p/q)^{1/3}+O(\mu)$,
$G(0)=(p/q)^{1/3}(1-e^2)^{1/2}$, $l(0)=n_l\pi$, and $g(0)=n_g\pi$,
where $p$ and $q$ are relatively prime positive integers,
$0<e<1$, and $n_l$ and $n_g$ are either $0$ or $1$, provided the solution
at $\mu=0$ does not collide with the orbit of the primary of mass $\mu$.
The existence holds for $\mu$ sufficiently small, the dependence
on $\mu$ is analytic, and there is only one such family as proved
in \cite{Arenstorf1}, \cite{Barrar1}. Given $p$, $q$, and $e$ there are four
possible choices for $n_l$ and $n_g$, but only two of these give rise to
distinct families. Two of the characteristic multipliers are equal to
$1\pm \sqrt{C(e,p,q)\mu}+O(\mu)$ in the limit $\mu\rightarrow 0$, where
$C(e,p,q)$ is given by \eqref{eqn-2-5}. The quantity $C(e,p,q)$ is analytic
at $e=0$ and $C(e,p,q) = O(e^{\abs{p-q}})$ as $e\rightarrow 0$, if
$p/q\neq 1/1$. The coefficient in front of $e^{\abs{p-q}}$ in the power 
series of $C(e,p,q)$ is given by \eqref{eqn-2-11}, \eqref{eqn-2-16},
and \eqref{eqn-2-17}. The coefficients for the two distinct families obtained
for given $p$, $q$, and $e$ sum to zero. 
\end{thm}

\begin{figure}
\begin{center}
\includegraphics[scale=0.35]{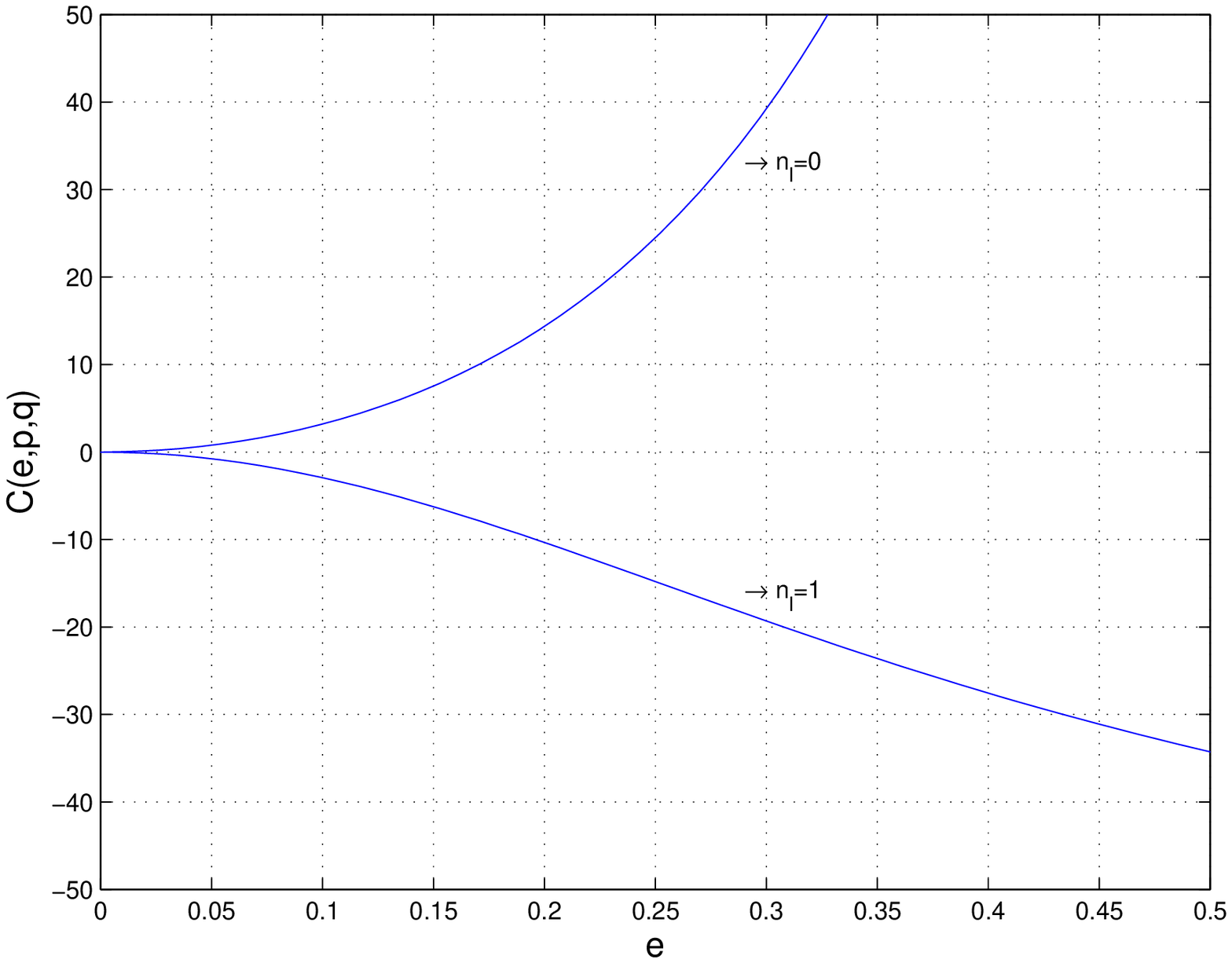}
\hspace{1cm}
\includegraphics[scale=0.35]{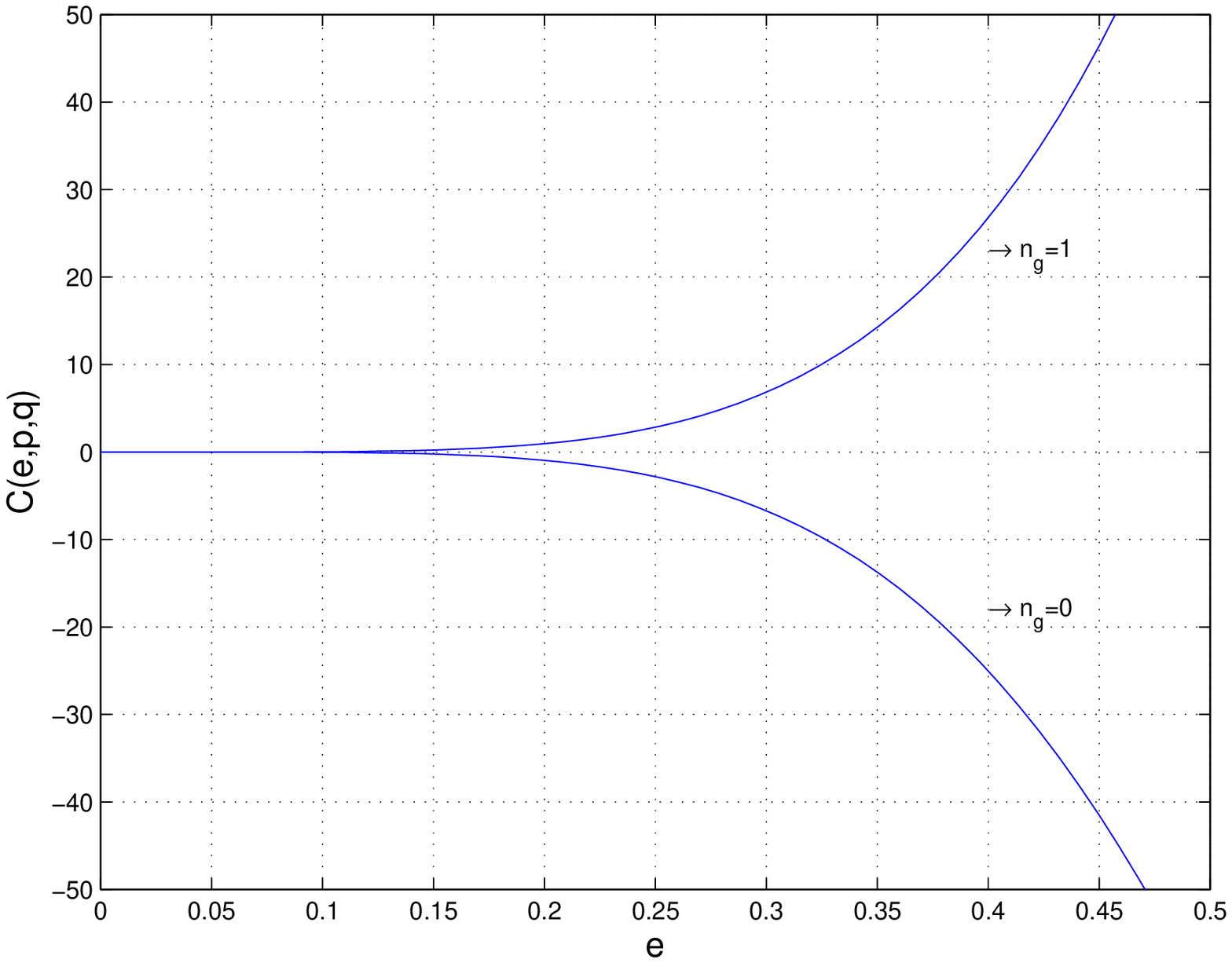}
\end{center}
\caption{The plots use $p/q=1/3$ and $p/q=2/7$, respectively}
\label{fig-3}
\end{figure}

The numerical evaluation of $C(e,p,q)$ defined by \eqref{eqn-2-5} can be carried
out with great accuracy at little expense. Figure \ref{fig-3} graphs
$C(e,p,q)$ against $e$ with two choices of $p/q$. In the graphs, $C(e,p,q)$
has opposite signs for the two distinct families. For $\mu > 0$ and
$\mu$ small, if the periodic solutions of one family are hyperbolic, the periodic
solutions of the other family are elliptic in either plot
of Figure \ref{fig-3}. This corresponds to the situation
depicted in Figure \ref{fig-2}. 

\section{Retrograde periodic motions}
The convention for representing retrograde motions in the restricted three body
problem depends upon the choice of variables. We will represent retrograde motions
using negative $L$ and $G$. The variation of $L,G,l,g$ with $t$ is again obtained 
from Hamilton's equations \eqref{eqn-2-1}.  The dependence of $r$ and $\theta$
on $L,G,l,g$ is given by \eqref{eqn-1-3} and \eqref{eqn-1-4} without a change.
For the choice of variables to be valid,
the condition $0<e<1$ must hold.
If the angles $E$, $\nu$, and $g$ are interpreted in the inertial frame,
say with $\mu=0$,
they must
all be measured in the counterclockwise sense. However, for retrograde motion, unlike
for direct motion, the angles $E$ and $\nu$ decrease with $t$. 

There exists a $\mu$-dependent family of retrograde periodic solutions of \eqref{eqn-2-1}
with initial conditions $L(0)=-(p/q)^{1/3}+O(\mu)$, $G(0)=-(p/q)^{1/3}(1-e^2)^{1/2}$,
$l(0) = n_l\pi$, and $g(0)=n_g\pi$, where $p$ and $q$ are relatively prime positive integers, $0<e<1$, and $n_l$ and $n_g$ are either $0$ or $1$. The dependence on
$\mu$ is analytic and the existence holds for $\mu>0$ and $\mu$ sufficiently small 
provided the solution at $\mu=0$ satisfies $\Delta_1 > 0$. Given $p$, $q$,
and $e$, there are four possible choices for $n_l$ and $n_g$. Of these only two
give rise to distinct families. If $p$ is odd, the two distinct families can be 
obtained with $n_g=0$ and $n_l$ either $0$ or $1$, and if $p$ is even with
$n_l=0$ and $n_g$ either $0$ or $1$.

Up to terms of order $\mu$, the dependence of periodic solutions of the family
determined by $p$, $q$, $e$, $n_l$, and $n_g$ is given by
\begin{align*}
L(t) &= -(p/q)^{1/3}+O(\mu)\\
G(t) &= -(p/q)^{1/3}(1-e^2)^{1/2} + O(\mu)\\
l(t) &= n_l\pi - qt/p + O(\mu)\\
g(t) &= n_g\pi -t + O(\mu)
\end{align*}
with the period being $T=2\pi p + O(\mu)$. To find the characteristic multipliers
of these periodic solutions, we choose the Poincar\'{e} section $g=n_g\pi$ and
represent it using the $L$-$l$ plane. As in Section 2, we may deduce that
\[\frac{\partial(L(T), l(T))}{\partial(L(0), l(0))} 
= \begin{pmatrix} 1 & 0\\
 -\frac{6\pi q^{4/3}}{p^{1/3}} & 1 \end{pmatrix}
 + \mu \begin{pmatrix} \star & -C(e,p,q)\frac{p^{1/3}}{6\pi q^{4/3}}\\
 \star & \star \end{pmatrix}
 + O(\mu^2).\]
The quantity $C(e,p,q)$ is given by
\begin{equation}
C(e,p,q) = -\frac{6\pi p^{5/3}}{q^{2/3}}\int_0^{2\pi p}
\Omega_{gg}\bigl(-(p/q)^{1/3},-(p/q)^{1/3}(1-e^2)^{1/2}, n_l\pi-qt/p, n_g\pi-t\bigr)
\, dt.
\label{eqn-3-1}
\end{equation}
This expression for $C(e,p,q)$ is similar in form to \eqref{eqn-2-4}, but 
the arguments of $\Omega_{gg}$ in \eqref{eqn-2-4} and \eqref{eqn-3-1} are
different. Use the change of variables from $t$ to $F$ given by 
$l=-qt/p+n_l\pi$, $l=E-e\sin E$, $dt = -(p/q)^{1/3}r dE$, and $qF=E$, to
get 
\begin{equation}
C(e,p,q) = -6\pi p^2 \int_0^{2\pi} \Bigl(\frac{r}{\Delta_1}\Bigr)_{\theta\theta}
+ \frac{\cos\theta}{r} \, dF,
\label{eqn-3-2}
\end{equation}
where $\Delta_1 = (1+r^2-2r\cos\theta)^{1/2}$. This expression for $C(e,p,q)$
looks identical to \eqref{eqn-2-6}, but the dependence of $\theta$ on 
$t$, and hence on $F$, is quite different. The dependence of $L,G,l,g$ on $t$
here is as indicated by the arguments of $\Omega_{gg}$ in \eqref{eqn-3-1}
and not as in \eqref{eqn-2-4}.  The variables $r$ and $\theta$ are obtained
from $L,G,l,g$ using \eqref{eqn-1-3} and \eqref{eqn-1-4} as before.
We now consider $C_1(e,p,q)$ and $C_2(e,p,q)$ which are defined as in
\eqref{eqn-2-7} but with the dependence of $r$ and $\theta$ on $F$ that is indicated
here. We have $C(e,p,q) = -6\pi p^2(C_1(e,p,q)+C_2(e,p,q))$.

That $C(e,p,q)$ is analytic in $e$ follows from the assumption $\Delta_1>0$ as in
Section 2. Here too $C(e,p,q)$ is analytic at $e=0$ if $p/q\neq 1/1$.
But the power series of $C(e,p,q)$ about $e=0$ looks quite different, as will be
shown now.   
 
Let $z=\exp(iF)$. 
Define $\beta$ by $e=2\beta/(1+\beta^2)$ as before.
We have 
\begin{align}
\exp(i\nu) &= z^q(1-\beta z^{-q})(1-\beta z^q)^{-1}\nonumber\\
\exp(-it) &= \exp(-i\pi n_lp/q)z^p\exp\Bigl(\frac{ep}{2q}(z^{-q}-z^q)\Bigr),
\label{eqn-3-3}
\end{align}
which has $-t$ in place of the $t$ in \eqref{eqn-2-8}.  From
$r = (p/q)^{1/3}(1-e\cos E)$, $\theta = \nu+n_g\pi -t$, and \eqref{eqn-3-3},
we have
\begin{align}
r &= (p/q)^{2/3} (1+\beta^2)^{-1}(1-\beta z^{-q})(1-\beta z^q)\nonumber\\
\exp(in\theta) &= \exp(-i\pi nn_lp/q)(-1)^{nn_g}z^{n(p+q)}(1-\beta z^{-q})^n
(1-\beta z^q)^{-n}\exp\bigl(\frac{enp}{2q}(z^{-q}-z^q)\Bigr).
\label{eqn-3-4}
\end{align}
The expression for $\exp(in\theta)$ in \eqref{eqn-3-4} has a factor $z^{n(p+q)}$
while the expression in \eqref{eqn-2-9} has the factor $z^{n(q-p)}$. A consequence
of the difference in the two expressions is that $C(e,p,q)$ for retrograde
motion is $O(e^{p+q})$,
and not $O(e^{\abs{p-q}})$, in the limit $e\rightarrow 0$ for $p/q\neq 1/1$, as will
be shown now.

The analysis of $C_2(e,p,q)$ is similar to that in Section 2 with the difference
that \eqref{eqn-3-4} and not \eqref{eqn-2-9}
must be used to express $\cos(\theta)$ in terms of $z$. We can conclude that
$C_2(e,p,q)=0$ if $q\neq 1$ and that $C_2(e,p,q) = O(e^{p+1})$ when $q=1$.
Besides,the sum of the coefficients of $e^{p+1}$ in the power series of $C(e,p,q)$
about $e=0$ for the two distinct families, obtained for given
$p$, $q$, and $e$, is always zero.
 
We consider $C_1(e,p,q)$ with the assumption $p<q$. The $p>q$ case is similar and leads
to the same conclusions. The integrand $\bigl(\frac{r}{\Delta_1}\bigr)_{\theta\theta}$
can again be written as a sum of terms as in \eqref{eqn-2-14} and the operator
$D_\alpha$ is again given by \eqref{eqn-2-13}. But the expression for the 
$n$th term in \eqref{eqn-2-14} is now different. Instead of \eqref{eqn-2-15}, we now
have
\begin{equation}
-n^2D_\alpha \alpha b_n(\alpha) \exp(in\theta) 
= -(-1)^{nn_g}\exp(-i\pi nn_l p/q) z^{n(p+q)} X_n(-D,D+n,D-n)(\alpha b_n(\alpha)),
\label{eqn-3-5}
\end{equation}
where $X_n(A,B,C) = (1+\beta^2)^A (1-\beta z^{-q})^B (1-\beta z^q)^C
\exp\bigl(\frac{enp}{2q}(z^{-q}-z^q)\bigr)$. The Laurent series of the 
$n$th term of \eqref{eqn-2-14}, which is displayed in \eqref{eqn-3-5},
has a possibly nonzero constant term only if $n=n_0q$, where $n_0$ is an integer,
and $n_0\neq 0$. That constant term is $\pm n^2$ times the coefficient of
$z^{n_0(p+q)q}$ in $X_n(-D,D+n,D-n)$, and therefore, if that term is expanded
about $e=0$ the lowest term with a possibly nonzero coefficient is 
$e^{\abs{n_0(p+q)}}$. Thus the lowest term with a possibly nonzero coefficient in
the power series of $C_1(e,p,q)$ about $e=0$ is $e^{p+q}$ and the only terms
\eqref{eqn-2-14} which contribute to its coefficient are obtained by setting
$n=\pm q$ in \eqref{eqn-3-5}. We conclude that $C_1(e,p,q) = O(e^{p+q})$
as $e\rightarrow 0$ for $p/q\neq 1/1$.
 
For given $p$, $q$, and $e$, the two distinct families are given by $n_g=0$
and $n_l$ either $0$ or $1$ if $p$ is odd, and by $n_l=0$ and
$n_g$ either $0$ or $1$ if $p$ is even. Inspection of \eqref{eqn-3-5} with
$n=\pm q$ leads to the conclusion that the sum of the coefficients of 
$e^{p+q}$ in the power series of $C_1(e,p,q)$ about $e=0$ for the two 
families is zero.

The theorem below is about retrograde periodic motions.

\begin{thm}
There exists a $\mu$-dependent family of periodic solutions of Hamilton's
equations \eqref{eqn-2-1} with the initial conditions 
$L(0)=-(p/q)^{1/3}+O(\mu)$, $G(0)=-(p/q)^{1/3}(1-e^2)^{1/2}$,
$l(0) = n_l \pi$, and $g(0)=n_g\pi$, where $p$ and $q$ are relatively prime positive
integers, $0<e<1$, and $n_l$ and $n_g$ are either $0$ and $1$, provided the
solution at $\mu=0$ does not collide with the orbit of the primary of mass $\mu$.
The existence holds for $\mu$ sufficiently small, the dependence on $\mu$ is
analytic, and there is only one such family. Given $p$, $q$, and $e$, the four
possible choices for $n_l$ and $n_g$ give rise to only two distinct families. 
Two of the characteristic multipliers are equal to 
$1\pm \sqrt{C(e,p,q)\mu}+O(\mu)$ in the limit $\mu\rightarrow 0$, where
$C(e,p,q)$ is given by \eqref{eqn-3-2}. The quantity $C(e,p,q)$ is analytic
at $e=0$ and $C(e,p,q) = O(e^{p+q})$ as $e\rightarrow 0$, if $p/q\neq 1/1$.
The sum of the coefficients of $e^{p+q}$ in the power series of $C(e,p,q)$
about $e=0$ for the two families is zero.
\end{thm}
  
\section{Periodic motions near collision}
In this section, we consider periodic motions near collision with the primary of mass
$1-\mu$. Consider the Hamiltonian $H$ defined by \eqref{eqn-1-0}. The generating
function $S = (-\mu+\xi^2-\nu^2)p_x+(2\xi\nu)p_y$ can be used to effect the
Levi-Civita transformation from the variables $p_x,p_y,x,y$ to the
variables $p_\xi, p_\nu,
\xi, \nu$. The Hamiltonian becomes
\[ H = \frac{p_\xi^2+p_\nu^2}{8(\xi^2+\nu^2)}+\frac{1}{2}(\nu p_\xi - \xi p_\nu)
+\frac{\mu}{2}\frac{(\xi p_\nu+\nu p_\xi)}{(\xi^2+\nu^2)}-\frac{1-\mu}{\xi^2+\nu^2}
-\frac{\mu}{\bigl((\xi^2-\nu^2-1)^2+4\xi^2\nu^2\bigr)^{1/2}}.\]
Suppose we are interested only in solutions of Hamilton's equations with $H=C$.
Consider $K=(\xi^2+\nu^2)(H-C)$ or
\[K = \frac{p_\xi^2+p_\nu^2}{8}+\frac{\xi^2+\nu^2}{2}(\nu p_\xi-\xi p_\nu-2C)-1
+\frac{\mu}{2}(\xi p_\nu+\nu p_\xi)+\mu -\frac{\mu(\xi^2+\nu^2)}
{\bigl((\xi^2-\nu^2-1)^2+4\xi^2\nu^2\bigr)^{1/2}}.\]
The solutions of Hamilton's equations formed using $K$ with $K=0$ correspond
to  solutions of Hamilton's equations of $H$ with $H=C$, but with time rescaled
from $t$ to $\tau$ such that $dt = (\xi^2+\nu^2)d\tau$. Hamilton's equations 
of
$K$ admit solutions with $K=0$ that pass through or close to $\xi=\nu=0$.
These can be interpreted as analytic continuations of the solutions of the equations
formed using $H$ through the singularity of $H$ at the location of the primary of
mass $1-\mu$. When $\mu=0$
\begin{equation}
K = \frac{p_\xi^2+p_\nu^2}{8}+\frac{\xi^2+\nu^2}{2}(\nu p_\xi-\xi p_\nu-2C)-1.
\label{eqn-4-1}
\end{equation}
The quantity $G=\xi p_\nu-\nu p_\xi$ is a first integral of the Hamilton's equations
of K shown in \eqref{eqn-4-1}. 
If a solution of those equations passes through $\xi=\nu=0$, then $G=0$.
In terms of $p_x, p_y, x, y$, $G=2(xp_y-yp_x)$.
Therefore the quantity denoted by $G$ in this section is twice the angular momentum.

The Hamilton-Jacobi equation of K shown in \eqref{eqn-4-1}was solved in \cite{Giacaglia1}
in the region of phase space 
with $G=0$ or $G\approx 0$. There is a gap in that solution 
when $G=0$. In this section, we discuss that gap and make  minor corrections to the solution given in \cite{Giacaglia1} for $G\neq 0$. We follow the geometric approach to the construction of action-angle variables described in \cite{Arnold1}. All references to
\cite{Arnold1} in this section are to the last chapter of that book.

The generating function $S=p_\xi (r\cos \theta)+p_\nu (r\sin\theta)$ can be used
to effect a change from the variables $p_\xi, p_\nu, \xi, \nu$ to the variables
$R, G, r, \theta$. The change to the polar variables $r,\theta$ and the corresponding
generalized momenta $R,G$ is well defined only if $r\neq 0$, or equivalently,
$\xi^2+\nu^2 \neq 0$. The Hamiltonian $K$ becomes
\[ K = \frac{1}{8}\Bigl(R^2+\frac{G^2}{r}\Bigr)+\frac{r^2}{2}(-G-2C)-1
+\frac{\mu}{2}(Rr\sin(2\theta)+G\cos(2\theta)) + \mu 
-\frac{\mu r^2}{(1+r^4-2r^2\cos(2\theta))^{1/2}}.\]
When $\mu=0$
\begin{equation}
K = \frac{1}{8}\Bigl(R^2+\frac{G^2}{r}\Bigr)+\frac{r^2}{2}(-G-2C)-1.
\label{eqn-4-2}
\end{equation}
We will find the action-angle variables of $K$ at $\mu=0$ using \eqref{eqn-4-2}.

The action-angle variables will be found in the region of $R,G,r,\theta$ space
or of $p_\xi, p_\nu, \xi, \nu$ space 
where the following conditions are satisfied:
\begin{equation}
G+2C < 0,\quad K+1 > 0, \quad (K+1)^2+\frac{G^2(G+2C)}{4} > 0,
\label{eqn-4-3}
\end{equation}
with $C<0$.
The necessity of these conditions will become clear shortly. We are interested
in solutions that are near a collision with the primary of mass $1-\mu$. Therefore
$G\approx 0$ in the region of interest. Further, the solutions of interest are bounded
and satisfy $K=0$. Therefore $C<0$ and $K+1\approx 1$ in the region of interest. It
is obvious that all points in the region of interest meet the conditions \eqref{eqn-4-3}.

\begin{figure}
\begin{center}
\begin{pspicture}(-3,-3)(3,3)
\pscircle[fillstyle=solid,fillcolor=gray]{2.7}
\pscircle[fillstyle=solid,fillcolor=white]{1.0}
\psline[arrowsize=3pt 2]{->}(0,0)(3.2,0)
\psline[arrowsize=3pt 2]{->}(0,0)(2.12,2.12)
\psarc{->}{0.5}{0}{45}
\uput[ur](0.55,0){$\theta$}
\uput[u](3.2,0){$r$}
\psarc[linewidth=2pt,arrowsize=3pt 2]{->}{1.85}{20}{380}
\pscurve[linewidth=2pt]{->}(1.0,0)(0.9,-0.1)(0.8,-0.2)(2.9,-0.2)(2.8,-0.1)(2.7,0)
\end{pspicture}
\end{center}
\caption{The cycles used to obtain action-angle variable are shown above.}
\label{fig-4}
\end{figure}
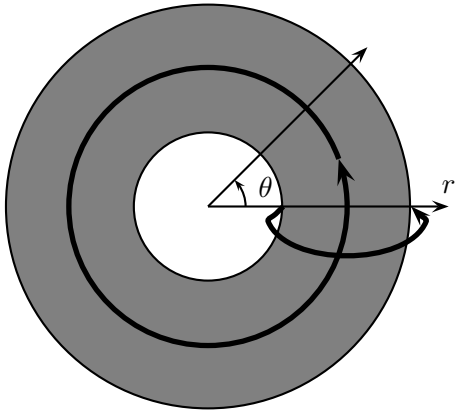

As a first step to the derivation of action-angle variables, we visualize the set
of points in $R,G,r,\theta$ for fixed values of $K$ and $G$ satisfying
\eqref{eqn-4-3} and $G\neq 0$. Using \eqref{eqn-4-2}, we deduce that
$r_{min} \leq r \leq r_{max}$ where
\begin{equation}
r_{min}^2+r_{max}^2 = \frac{2(K+1)}{(-G-2C)}, \quad
r_{min}^2 r_{max}^2 = \frac{G^2}{4(-G-2C)}.
\label{eqn-4-4}
\end{equation}
The last of the conditions \eqref{eqn-4-3} ensures $0\leq  r_{min} < r_{max}$.
The condition $G\neq 0$ implies $0< r_{min}$. For later use, we introduce variables
$a$ and $e$ defined by
\begin{equation}
r_{min}^2 = a(1-e), \quad r_{max}^2 = a(1+e).
\label{eqn-4-5}
\end{equation}
For $r$ within the allowed range, any value of $\theta$ is allowed in the set.
Thus the projection of the set looks like the annulus in Figure \ref{fig-4}.
The variable $G$ is of course fixed. For a given point in the annulus in the
$r$-$\theta$ plane, $R$ can have two values of equal magnitude but opposite sign.
On the boundary of the annulus, $R=0$. Therefore the set of points is a torus.

The set of points with $G=0$ and a fixed value of $K$, such that the conditions
\eqref{eqn-4-3} are met, includes points with $\xi=\nu=0$.
The change to polar variables is invalid at those points and the $R,G,r,\theta$
space
cannot be used to visualize that set. By \eqref{eqn-4-1}, the set of points
in the $p_\xi, p_\nu, \xi, \nu$ space is given by
the conditions
\[\xi p_\nu - \nu p_\xi=0, \quad
\xi^2 + \nu^2 \leq (K+1)/(-C), \quad
p_\xi^2+p_\nu^2 = 8(K+1+C(\xi^2+\nu^2)).\]
This set of points is also homeomorphic to a torus. It is easy to check that
the differential $1$-forms $dG$ and $dK$ are linearly independent on this torus.

We return to Figure \ref{fig-4} and the assumption $G\neq 0$ and find the 
action-angle variables. As explained in \cite{Arnold1}, the action variables
are functions of $K$ and $G$ obtained as
$\frac{1}{2\pi}\int R\,dr + G\, d\theta$, with the integral taken
over two linearly independent cycles on the torus. Our choice of cycles is shown
in Figure \ref{fig-4}. The cycle with fixed $r$ and $\theta$ varying from
$0$ to $2\pi$ gives $G$ as an action variable. The other cycle increases from
$r_{min}$ to $r_{max}$ with $R\geq 0$ and then decreases to $r_{min}$ with
$R\leq 0$. The variable $\theta$ is fixed on this cycle. The action variable
$L^{\star}$ given by this cycle is 
\[ L^{\star} = \frac{1}{\pi}\int_{r_{min}}^{r_{max}}
\Bigl(8K+8-\frac{G^2}{r^2}+4 r^2 (G+2C)\Bigr)^{1/2} dr.\]
Substitute $u = r^2$ to get
\[ L^{\star} = \frac{1}{\pi} (-G-2C)^{1/2} \int_{a(1-e)}^{a(1+e)}
\frac{\bigl((a(1+e)-u)(u-a(1-e))\bigr)^{1/2}}{u} du. \]
Evaluate the integral and use \eqref{eqn-4-4} and \eqref{eqn-4-5} to get
\begin{equation}
L^{\star} = \frac{K+1}{(-G-2C)^{1/2}}-\frac{\abs{G}}{2}.
\label{eqn-4-6}
\end{equation}
Instead of $L^{\star}$, we use $L$ defined by
\begin{equation}
L = \frac{K+1}{(-G-2C)^{1/2}}
\label{eqn-4-7}
\end{equation}
as the other action variable. 

The angle variables parametrize the surface of the torus obtained in $R,G,r,\theta$
space by fixing the action variables $L$ and $G$. As proved in \cite{Arnold1}, the
angle variables can be derived from the following generating function:
\[S(r,\theta, L, G) = \int_{r_{min},\,0}^{r,\,\theta} G\,d\theta + R\,dr.\]
In this expression, $R$ must be expressed in terms of $L$, $G$, and $r$
using \eqref{eqn-4-2} and \eqref{eqn-4-7}. 
From this expression, it might seem that $S(r,\theta,L,G)$ is a function on the
annulus  shown in Figure \ref{fig-4} for fixed $L$ and $G$. It is actually a function
on the torus as $R$ can be either positive or negative.
\[S(r,\theta,L,G) = G\theta \pm \int_{r_{min}}^r \Biggl(8L(-G-2C)^{1/2}
-\frac{G^2}{r^2} + 4 r^2(G+2C)\Biggr)^{1/2} dr,\]
where the sign is $+$ if $R\geq 0$ and $-$ otherwise. The generating function
$S(r,\theta,L,G)$ is a multiple valued function on the surface of the torus
as there are many non-homotopic paths from the base point $R=0, r=r_{min}, \theta=0$
to any point on the torus. In the calculations below, we consider the paths for
which both $r$ and $\theta$ increase monotonically.

The value of one of the angle variables on a point on the torus is equal to
the value of
$\frac{\partial S}{\partial L}$ at that point on the torus as proved in \cite{Arnold1}.
\[\frac{\partial S}{\partial L} =
\pm 4(-G-2C)^{1/2} \int_{r_{min}}^{r} \Biggl(8L(-G-2C)^{1/2}-\frac{G^2}{r^2}
+ 4r^2(G+2C)\Biggr)^{-1/2} dr,\]
where the sign is $+$ if $R\geq 0$ at the point on the torus and $-$ otherwise.
Change variable to $u=r^2$ to get
\[ \frac{\partial S}{\partial L} =
\pm \int_{a(1-e)}^{r^2} \frac{du}{\bigl((a(1+e)-u)(u-a(1-e))\bigr)^{1/2}}\,.\]
Change variable to $l$ defined by $u=a(1-e \cos l)$ with $0\leq l \leq \pi$
to get $\frac{\partial S}{\partial L} = \pm l$. If we adopt the convention 
$0\leq l \leq \pi$ if $R\geq 0$ and $-\pi< l < 0$ if $R< 0$, we get
\begin{equation}
\frac{\partial S}{\partial L} = l 
\label{eqn-4-8}
\end{equation}
where $l$ is an angular variable that is measured modulo $2\pi$.
 
The value of the other angle variable at a point on the torus 
is equal to the value of $\frac{\partial S}{\partial G}$ at that point on the
torus.
\[\frac{\partial S}{\partial G} = \theta 
\pm \frac{1}{2}\int_{r_{min}}^r \Biggl(-4L(-G-2C)^{-1/2}-\frac{2G}{r^2}+4r^2\Biggr)
\Biggl(8L(-G-2C)^{1/2}-\frac{G^2}{r^2}+4r^2(G+2C)\Biggr)^{-1/2} dr,\]
where the sign is $+$ if $R\geq 0$ at the point on the torus and $-$ otherwise.
Change variable to $u=r^2$ to get
\[\frac{\partial S}{\partial G} = \theta
\pm \frac{1}{4(-G-2C)^{1/2}}\int_{a(1-e)}^{r^2}
\frac{\bigl(-2L(-G-2C)^{-1/2}-G/u+2u\bigl)}
{\bigl((a(1+e)-u)(u-a(1-e)\bigr)^{1/2}}\,du.\]
Change variable to $l$ defined by $u=a(1-e\cos l)$ with the same convention
as in \eqref{eqn-4-8}.
By \eqref{eqn-4-4}, \eqref{eqn-4-5}, and \eqref{eqn-4-7},
\[-2L(-G-2C)^{1/2}-G/u+2u = -\frac{G}{a(1-e\cos l)}-2ae\cos l.\]
Denote $\frac{\partial S}{\partial G}$ by $g$ to get 
\begin{equation}
g = \theta - \frac{G}{4L} \int_0^l \frac{dl}{(1-e\cos l)} - \frac{(L^2-G^2/4)^{1/2}}
{2(-G-2C)}\sin l,
\label{eqn-4-9}
\end{equation}
where $e=\bigl(1-\frac{G^2}{4L^2}\bigr)^{1/2}$.

If $L^{\star}$, defined by \eqref{eqn-4-6}, and $G$ are used as the action variables,
the corresponding angle variables are $l$ and $g\pm l/2$
--- the sign is $+$ if $G>0$ and $-$ if $G<0$.
It is easily verified that
along each of the two cycles depicted in Figure \ref{fig-4} one of these angle
variables increases by $2\pi$ while the other doesn't change. This must be the
case as proved in \cite{Arnold1} and the verification is a check of the correctness
of our derivation of action-angle variables. 
 
The quantities denoted by $G$ and $g$ in \cite{Giacaglia1} are twice and half
of the quantities denoted by $G$ and $g$ in this section, respectively. 
The factor $2$ that appears in the coefficient
of $\sin l$ in (19) of \cite{Giacaglia1} must be moved to the denominator.
Further, in (20) of \cite{Giacaglia1} $(1-e^2)^{1/2}$ must be replaced by
$G/L$. Although $\abs{G/L} = (1-e^2)^{1/2}$ in the notation of 
\cite{Giacaglia1}, the sign of $G$ is significant.

When $G=0$ the change to polar variables is not valid as we have shown in this section.
Therefore the derivation of action-angle variables using a change to polar
variables given in \cite{Giacaglia1} is also not valid when $G=0$.

The existence of $\mu$-dependent families of periodic solutions of Hamilton's
equations of $K$ with the initial conditions $L(0) = (p/q)^{1/3}+O(\mu)$,
$G(0)=-2C-(q/p)^{2/3}+O(\mu)$, $l(0) = 0$, and $g(0)=0$ was proved
in \cite{Schmidt1}. The solutions in these families depend analytically on $\mu$
and $K=0$ at all points along the solutions. 
The methods of Sections 2 and 3 can be used to carry out linear stability analyses
of these periodic solutions. 
However, a new derivation of 
action-angle variables is necessary for the crucial $G=0$ case.

%\cite{Arenstorf1} \cite{Barrar1} \cite{Schmidt1} \cite{MH1} \cite{SM1} \cite{Goursat1}
%\cite{LS1} \cite{BC1} \cite{Giacaglia1} \cite{Birkhoff1} \cite{Arnold1}.
\bibliography{references}
\bibliographystyle{plain}

\end{document}